\documentclass[11pt, leqno]{article}
\usepackage{graphicx}
\usepackage{amsmath,amssymb,amscd,amsthm, latexsym}

\newtheorem{theorem}{Theorem}[section]

\newtheorem{cor}[theorem]{Corollary}

\newtheorem{conjecture}[theorem]{Conjecture}

\theoremstyle{definition}

\theoremstyle{remark}

\newcommand{\address}{Department of Mathematical Sciences,
        Korea Advanced Institute of Science and Technology,
        Yusong-gu, Daejon 305--701, Republic of Korea\par
        {\it E-mail address}: {\tt jinkim@math.kaist.ac.kr}
        }

\numberwithin{equation}{section}

\begin{document}

\title%%[On positively curved 4-manifolds with $S^1$-symmetry]
{\bf On positively curved 4-manifolds with $S^1$-symmetry}

\author{Jin Hong Kim}

%%\date{\today} %%: December 24, 2005

%%\subjclass{Primary 53C15}

%%\keywords{positively curved 4-manifolds, effective circle action,
%%diffeomorphism classification}

%%\thanks{This work was supported by KOSEF.}

\maketitle

% ----------------------------------------------------------------
\begin{abstract}
It is well-known by the work of Hsiang and Kleiner that
every closed oriented positively curved 4-dimensional manifold with
an effective isometric $S^1$-action is homeomorphic to $S^4$ or
${\bf CP}^2$. As stated, it is a topological classification. The primary
goal of this paper is to show that it is indeed a diffeomorphism
classification for such 4-dimensional manifolds. The proof of this diffeomorphism classification also shows an even stronger statement that every positively curved simply connected $4$-manifold with an isometric
circle action admits another smooth circle action which extends to a $2$-dimensional torus action and is equivariantly diffeomorphic to a linear action on $S^4$ or ${\bf CP}^2$.
The main strategy is to analyze all possible topological configurations of effective
circle actions on simply connected 4-manifolds by using the
so-called replacement trick of Pao.
\end{abstract}

%%\setlength{\baselineskip}{15pt}

% ----------------------------------------------------------------
\section{Introduction and Main Results} \label{sec1}

A closed Riemannian $n$-dimensional manifold with positive sectional
curvature everywhere is called a \emph{positively curved
$n$-manifold}. These manifolds are very rare and, moreover, satisfy
very special topological properties. In particular, a well-known
theorem of Synge says that every even dimensional orientable closed
positively curved manifold is simply connected, while every such odd
dimensional manifold has at most finite first fundamental group (see
\cite{DoC}). Hence it is immediate that in two dimensions two-sphere
is only an orientable closed positively curved manifold. In three
dimensions, Hamilton has classified all orientable closed positively
curved manifolds by using the techniques of Ricci flow \cite{H82}.
It turns out that they are all diffeomorphic to space forms. On the
other hand, in four dimensions the classification of orientable
closed positively curved manifolds is still incomplete. Much
interestingly, all the well-known orientable closed positively
curved 4-manifolds always admit a positively curved metric with a
circle symmetry.

In their paper \cite{HK89}, Hsiang and Kleiner investigated the
question which orientable closed positively curved 4-manifolds admit
a positively curved metric with an effective isometric $S^1$-action.
They showed that if $M$ is a compact oriented positively curved
4-manifold with an effective isometric $S^1$-action, then $M$ is
homeomorphic to $S^4$ or ${\bf CP}^2$. This classification is
remarkable and sparked off many current research activities in this
field. In particular, see the works of Grove--Searle and B. Wilking
in \cite{GS94} and \cite{W03, W05} respectively. However, it is a
topological classification rather than a diffeomorphism one. In
fact, Hsiang and Kleiner (and probably many others) asked in the
same paper whether or not the following conjecture was true:

\begin{conjecture} \label{conj1.1}
An orientable closed positively curved 4-manifold $M$ with an
effective isometric $S^1$-action is diffeomorphic to $S^4$ or ${\bf
CP}^2$.
\end{conjecture}

The main aim of this paper is to affirmatively prove their Conjecture
\ref{conj1.1} (or Conjecture 1 in \cite{HK89}). What we have done to
prove the conjecture is just to add some more remarks to the
homeomorphism classification of Hsiang and Kleiner.

Throughout this paper, all $S^1$-actions are assumed to be effective,
unless stated otherwise. In order to explain the main ideas, we let
$F(S^1, M)$ be the fixed point set of such an $S^1$-action on $M$.
It is well known as in \cite{Kob} that the Euler characteristic of
$F(S^1, M)$ is equal to that of $M$ in the presence of an effective
$S^1$-action. Thus the Euler characteristic of the fixed point set
$F(S^1, M)$ is greater than or equal to 2. One of the main steps in
the paper of Hsiang and Kleiner is to prove that the Euler
characteristic of the fixed point set $F(S^1, M)$ is at most three.
Hence, since each component of $F(S^1, M)$ is a totally geodesic
submanifold of $M$, $F(S^1, M)$ has the following four possibilities
only:
\begin{enumerate}

\item[(1)] One 2-sphere.

\item[(2)] The disjoint union of one 2-sphere and one isolated fixed point.

\item[(3)] Two isolated fixed points.

\item[(4)] Three isolated fixed points.

\end{enumerate}

Hsiang and Kleiner obtained the above information about the fixed
point set $F(S^1, M)$ by essentially using the existence of a
positively curved metric with an effective isometric $S^1$-action.
But it seems to need more serious consideration on the topological
properties of effective circle actions on simply connected
4-manifolds. Thus our starting point of this paper is to investigate
all possible topological configurations of effective circle actions
on simply connected 4-manifolds. However, the existence of a
positively curved metric with an effective isometric $S^1$-action is
also used crucially throughout. Moreover, we need to use the recent
resolution of the Poincar\' e conjecture by Perelman as in
\cite{Per1} and \cite{Per2}. (See also \cite{M-T}, \cite{C-Z, C-Z2}
and \cite{K-L}.) So our proof is differential-topological in nature.
We remark that it is inspired by the work \cite{B02} of S.
Baldridge.

In case that the fixed point set $F(S^1, M)$ is either one 2-sphere
or the union of one 2-sphere and one isolated fixed point, the
effective isometric $S^1$-action on $M$ has the fixed point set
whose codimension is 2. It then follows from Theorem 1.2 of Grove
and Searle in \cite{GS94} that $M$ is diffeomorphic to $S^4$ or
${\bf CP}^2$. So it is enough to consider the remaining two cases:
$F(S^1, M)$ consists of either two isolated fixed points or three
isolated fixed points. In the present paper, we prove that even in
this case $M$ is diffeomorphic to $S^4$ or ${\bf CP}^2$.

One of the new ingredients of this paper that is not present in the
paper of Hsiang and Kleiner is to use the classification results of
circle actions on simply connected 4-manifolds by R. Fintushel in
\cite{F77} and \cite{F78}. According to Fintushel, smooth
$S^1$-actions on simply connected 4-manifolds can be classified in
terms of their legally weighted orbit spaces. Applying this
classification to our situation, we will have at least four (resp.
three) possibilities for legally weighted orbit spaces, if the fixed
point set is three isolated fixed points (resp. two isolated fixed
points). By some case-by-case analysis we can finally show the
following result.

\begin{theorem} \label{thm1.1}
An orientable closed positively curved 4-manifold $M$ with an
effective isometric $S^1$-action is diffeomorphic to $S^4$ or ${\bf
CP}^2$.
\end{theorem}

This completes the classification of closed oriented positively
curved manifolds of dimension 4 with effective isometric
$S^1$-actions, up to diffeomorphism. Moreover, as the reader of this paper also pointed out, the proof of Theorem
\ref{thm1.1} in Section \ref{sec4} actually provides the following even stronger result.

\begin{theorem} \label{thm1.2}
A positively curved simply connected $4$-manifold $M$ with an isometric
circle action admits another smooth circle action which extends to a $2$-dimensional torus action and is equivariantly diffeomorphic to a linear action on $S^4$ or ${\bf CP}^2$.
\end{theorem}

When a fixed point set contains a $2$-dimensional component or there is no weight circle in the orbit space, Theorem \ref{thm1.2} can be further strengthened to the statement that an isometric circle action on a positively curved simply connected closed $4$-manifold is equivariantly diffeomorphic to a linear action on $S^4$ or ${\bf CP}^2$ (see also Theorem 2.8 of \cite{GS97} and Conjecture 3.1 in \cite{Ga-Gr}). So it would be interesting to know whether or not the strengthened statement still holds even in the presence of a weight circle in the orbit space.

We finally remark that the technique of this paper resolves a version of Theorem \ref{thm1.2} for nonnegatively curved simply connected 4-manifolds with an isometric $S^1$-action as well as the following conjecture in \cite{HK89} without any significant efforts (e.g., see \cite{K-Lee06}):

\begin{conjecture}
An orientable closed simply connected nonnegatively curved
4-manifold $M$ with an effective isometric $S^1$-action is
diffeomorphic to either $S^4$, ${\bf CP}^2$, ${\bf CP}^2\# \pm {\bf
CP}^2$, or $S^2\times S^2$.
\end{conjecture}

We organized this paper as follows. In Section 2, we review the main
ingredients in the classification of circle actions on simply
connected 4-manifolds done by Fintushel. This will be largely taken
from two papers \cite{F77} and \cite{F78}. In Section 3, we explain
the so-called \lq\lq replacement trick" by P. Pao in \cite{P78}. In
Section 4, we give detailed proofs of Theorems \ref{thm1.1} and \ref{thm1.2} by
considering the various topological configurations of the isolated
fixed points.

\section{Classification of circle actions on simply connected 4-manifolds} \label{sec2}

The aim of this section is to review the main ingredients in the
classification of circle actions on simply connected 4-manifolds by
Fintushel in \cite{F77} and \cite{F78}. These together with the
replacement trick by P. Pao in \cite{P78} will be our important
tools to prove Theorems \ref{thm1.1} and \ref{thm1.2}.

Throughout this section, let $M$ be a simply connected oriented
4-manifold with an effective $S^1$-action. For any subset $N$ of
$M$, let $N^\ast$ denote the image of $N$ in the orbit space and let
$p: M\to M^\ast$ be the orbit map. Also if we are given a set
$N^\ast$ in $M^\ast$, $N$ means the preimage of $N^\ast$ under the
orbit map $p$. Finally let $F$ be the fixed point set of $M$, $E$
the union of the exceptional orbits, and $P$ the union of the
principal orbits. Since $M$ is simply connected, $M^\ast$ is also
simply connected. Thus each component of $\partial M^\ast$ is a
2-sphere by the Poincar\' e duality and is a subset of $F^\ast$ of
the fixed point set. Then we can assign the orbital data to $M^\ast$
as follows:

\begin{enumerate}

\item[(1)]
For each boundary component $N^\ast$ of $M^\ast$, choose a regular
neighborhood $N^\ast \times [0,1]$. The restriction of the orbit map
$p$ to $N^\ast \times \{ 1 \}$ is a principal $S^1$-bundle over
$N^\ast \times \{ 1 \}$. Thus we assign to each boundary component
$N^\ast$ of $M^\ast$ oriented by the outward normal to $N^\ast
\times \{ 1 \}$ the Euler number of this bundle, and $N^\ast$ is
called a \emph{weighted 2-sphere}.

\item[(2)]
For each $x^\ast$ in $F^\ast$ minus the union of $\partial M^\ast$
and the closure of $E^\ast$, let $B^\ast$ be a 3-disk neighborhood
of $x^\ast$ in the interior of $M^\ast$. The restriction of $p$ to
the boundary $\partial B^\ast\subset P^\ast$ is a principal
$S^1$-bundle over $\partial B^\ast$ whose total space is a 3-sphere.
With the orientation of $\partial B^\ast$ by the outward normal of
$\partial B^\ast$, we assign to $x^\ast$ the Euler number $\pm 1$ of
this bundle.

\item[(3)]
For each simple closed curve $C^\ast$ in $E^\ast \cup F^\ast$, we
assign to each component of $E^\ast$ in $C^\ast$ the Seifert
invariant, once we fix an orientation in $C^\ast$. The simple closed
curve obtained in this way will be called \emph{multiply-weighted}.
If a weighted circle does not contain any fixed points, it is called
\emph{simply-weighted}. But if $M$ is simply connected, then
simply-weighted circles in the orbit space $M^\ast$ do not occur.

\item[(4)]
For each arc $A^\ast$ which is a component of $E^\ast\cup F^\ast$,
we orient $A^\ast$ and assign the Seifert invariant as in (3). In
this case we call $A^\ast$ a weighted arc and write the weight
system as
\[
\left[ b'; (\alpha_1, \beta_1), \ldots, (\alpha_n, \beta_n);
b''\right]
\]
\end{enumerate}

To each component of (1), (2), (3), and (4), their associated Euler
number will be called the \emph{index}. Thus the weighted arc
\[
\left[ b'; (\alpha_1, \beta_1), \ldots, (\alpha_n, \beta_n);
b''\right]
\]
has index $b''-b'=\pm 1$ or $0$, an isolated fixed point has index
equal to $\pm 1$, and a weighted circle has index equal to $0$.
Moreover, it is important to note that the sum of the indices of all
the components (1), (2), (3), and (4) is always zero. The main
result of \cite{F77} and \cite{F78} says that each simply connected
4-manifold $M$ with an effective $S^1$-action corresponds to a
\emph{legally weighted 3-manifold} $M^\ast$ which is an oriented
compact connected 3-manifold equipped with the following data:

\begin{enumerate}
\item[(1)]
A finite collection of weighted arcs and circles in the interior of
$M^\ast$ as above.

\item[(2)]
A finite collection of distinguished points in the interior of
$M^\ast$.

\item[(3)]
Let $S^\ast$ be the union of $\partial M^\ast$, the points of (2),
and the weight arcs of (1). A class $\chi$ is an element in
$H_1(M^\ast, S^\ast)$ such that if
\[
\partial \chi=(\zeta_1, \ldots, \zeta_m)\in H_0(S^\ast; {\bf Z})
\]
then each isolated fixed point has $\zeta_i=\pm 1$ and each weighted
arc
\[
\left[ b'; (\alpha_1, \beta_1), \ldots, (\alpha_n, \beta_n);
b''\right]
\]
has $\zeta_i=b''-b'=\pm 1$ or $0$.
\end{enumerate}

With these preliminaries in place, the following theorem (Theorem
7.1 in \cite{F77}) will play a crucial role in this paper.

\begin{theorem} \label{thm2.1}
Let $S^1$ act on a simply connected 4-manifold $M$. The action of
$S^1$ extends to an action of $T^2=S^1\times S^1$ if and only if
\begin{enumerate}
\item[\rm (1)]
$M^\ast$ is not a counterexample to the 3-dimensional Poincar\' e
conjecture and

\item[\rm (2)]
if $M^\ast$ contains a weighted circle $C^\ast$ then $M^\ast$ is
homeomorphic to $S^3$, $C^\ast= E^\ast\cup F^\ast$, and $C^\ast$ is
unknotted in $M^\ast$.
\end{enumerate}
\end{theorem}

In particular, this theorem can be used to show that if the two
conditions of Theorem \ref{thm2.1} are satisfied then a simply
connected 4-manifold with an effective $S^1$-action can be built
from $S^4$ with the circle action with two isolated fixed points and
equivariantly connect summing $S^4$'s and ${\bf CP}^2$'s repeatedly.
To be more precise, Fintushel proved the following theorem.

\begin{theorem} \label{thm2.2}
Let $S^1$ act on a simply connected 4-manifold $M$, and suppose that
the orbit space that is homotopy equivalent to $S^3$ is not a
counterexample to the 3-dimensional Poincar\' e conjecture. Then $M$
is equivariantly homeomorphic to a connected sum of copies of $S^4$,
$\pm{\bf CP}^2$, and $S^2\times S^2$.
\end{theorem}

It is now true by the celebrated works \cite{Per1} and
\cite{Per2} of Perelman that there is no counterexample to the
Poincar\' e conjecture. Hence we can state the following

\begin{cor} \label{cor2.1}
Let $S^1$ act on a simply connected 4-manifold $M$. Then $M$ is
equivariantly homeomorphic to a connected sum of copies of $S^4$,
$\pm {\bf CP}^2$, and $S^2\times S^2$.
\end{cor}

\section{Replacement Trick} \label{sec3}

\begin{figure}
\begin{center}
\includegraphics[height=60mm]{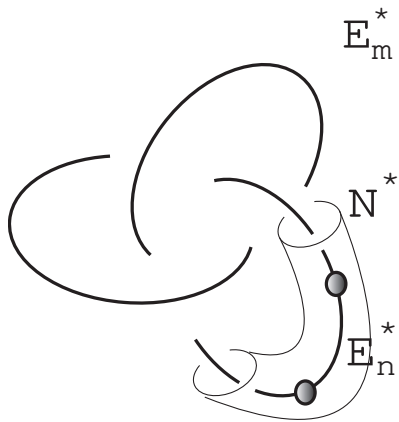}
\end{center}
\caption{A configuration for the replacement trick} \label{fig1}
\end{figure}

In this paper we also need the \lq\lq replacement trick" of P. Pao.
This trick makes the given 4-manifold admit many different
$S^1$-actions by replacing a weighted arc with a simpler one. To be
precise, if the weighted circle $C^\ast$ contains exactly two
isolated fixed points then $C^\ast$ consists of two exceptional
orbit types $E^\ast_n$ and $E^\ast_m$ and two fixed points. That is,
$C^\ast$ is just
\[
C^\ast= E^\ast_n\cup E^\ast_m\cup \{ \text{two isolated fixed
points} \}.
\]
Now choose a tubular neighborhood $N^\ast$ of the union of
$E^\ast_n$ and two fixed points such that $N^\ast$ intersect the
circle $C^\ast$ at an open arc containing two fixed points and
$E^\ast_n$. See Figure \ref{fig1}.

Let $N$ be the preimage of $N^\ast$ under the orbit map. Consider
the product $D^2\times S^2$ of a 2-disk and a 2-sphere. Then
represent $D^2\times S^2$ by $(r, \theta; \rho, \phi, z)$ where
$(r,\theta)$ are the polar coordinates of $D^2$ and $(\rho, \phi,
z)$ are the cylindrical coordinates of $S^2$. Let $S^1$ act on
$D^2\times S^2$ by
\[
(r, \theta; \rho, \phi, z)\mapsto (r, \theta+m\psi; \rho,
\phi+n\psi, z).
\]
Note that the orbit space of $D^2\times S^2$ under the above action
of $S^1$ is identical to $N^\ast$. On the other hand, there are many
other actions of $S^1$ on $D^2\times S^2$. For example, if $n'\equiv
n$ mod $m$ then we let $S^1$ act on $D^2\times S^2$ by
\[
(r, \theta; \rho, \phi, z)\mapsto (r, \theta+m\psi; \rho,
\phi+n'\psi, z).
\]
Now remove $N$ from $M$ and change the $S^1$-action on $N$. Then sew
$N$ back to $M$ by an equivariant diffeomorphism. If we require
$(n-n')/m$ to be even in addition, then this surgery operation will
not change the manifold $M$, but change the orbit structure on
$C^\ast$ by the union of two exceptional orbits $E^\ast_{n'}$ and
$E^\ast_m$ and two isolated fixed points. This reduction can be
repeated with keeping decrease the pair $(m,n)$ of integers to
$(1,1)$ or $(1,0)$. Hence we have the following theorem of P. Pao in
\cite{P78} (or Proposition 13.1 in \cite{F78}):

\begin{theorem} \label{thm3.1}
Let $X$ be a closed oriented 4-manifold with an effective
$S^1$-action whose weighted orbit space contains a weighted circle
$C^\ast$ with exactly two isolated fixed points. Then $M$ admits a a
different $S^1$-action whose weight space is $M^\ast$ with $C^\ast$
replaced by two isolated fixed points or $M^\ast$ minus the interior
of 3-ball with $C^\ast$ removed.
\end{theorem}

This implies that in the second case the orbit space for the new
$S^1$-action will have one more boundary component instead of the
weight circle $C^\ast$.

\section{Proofs of Theorems \ref{thm1.1} and \ref{thm1.2}} \label{sec4}

The goal of this section is to prove Theorems \ref{thm1.1} and \ref{thm1.2}. As
remarked in Section 1, it suffices to prove the theorems only for the
following two cases: two isolated fixed points or three isolated
fixed points. If the manifold $M$ has two (resp. three) isolated
fixed points, then $M$ is equivariantly homeomorphic to $S^4$ (resp.
${\bf CP}^2$) by Theorem \ref{thm2.2}.

\begin{figure}
\begin{center}
\includegraphics[height=20mm]{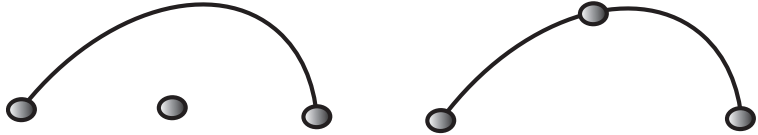}
\end{center}
\caption{Some configurations of $E^\ast\cup F^\ast$ for three
isolated fixed points}\label{fig2}
\end{figure}

Assume first that the orbit space $M^\ast$ does not have any
weighted circles, regardless of the number of isolated fixed points.
See Figure \ref{fig2} for some possible configurations of
$E^\ast\cup F^\ast$.

Since there is no counterexample to the Poincar\' e conjecture by
Perelman, it follows from Theorem \ref{thm2.1} that the $S^1$-action
lifts to an action of $T^2=S^1\times S^1$. Hence the simply
connected positively curved 4-manifold $M$ which is equivariantly
homeomorphic to $S^4$ or ${\bf CP}^2$ is actually a toric manifold.
But then the well-known classification result (or argument) of toric
manifolds immediately tells us that $M$ is in fact equivariantly
diffeomorphic to $S^4$ or ${\bf CP}^2$, depending on the number of
isolated fixed points. Or, to be more precise, you may adapt an argument
appearing at the end of this section which shows that a simply
connected manifold with an action of $T^2$ equivariantly
homeomorphic to ${\bf CP}^2$ is actually equivariantly diffeomorphic
to ${\bf CP}^2$.

Next we assume that the orbit space $M^\ast$ has a weighted circle.
Then we need to use the replacement trick of Pao (Theorem
\ref{thm3.1}). To do so, first consider the case that there are two
isolated fixed points. In this case $E^\ast\cup F^\ast$ consists of
one weighted circle with two isolated fixed points lying on it. Note
that there is no case that $F(S^1, M)$ is the disjoint union of a
simply weighted circle and two isolated fixed points, since $M$ is
simply connected (see \cite{F77} or Section 2). Now applying the
replacement trick to the manifold $M$ and its orbit space $M^\ast$,
we obtain either a new manifold $M'$ equivariantly diffeomorphic to
$M$ whose $E^\ast\cup F^\ast$ consists of two isolated fixed points
or a new manifold $M'$ whose fixed point set is a 2-sphere. In the
former case $M'$ is equivariantly diffeomorphic to $S^4$ as above.
On the other hand, in the latter case there exists an isometric
effective $S^1$-action on $M'$ with a positively curved metric whose
fixed point is of codimension 2. Thus again by Theorem 2.1 of Grove
and Searle in \cite{GS94} $M'$ (or $M$) should be equivariantly
diffeomorphic to $S^4$. In fact, for either case there exists an
alternative argument: we can apply Theorem \ref{thm2.1} to show that
$M'$ is equivariantly diffeomorphic to $S^4$, since anyway there is
no weighted circle by construction. This finishes the proof of the
case of two isolated fixed points.

\begin{figure}
\begin{center}
\includegraphics[height=60mm]{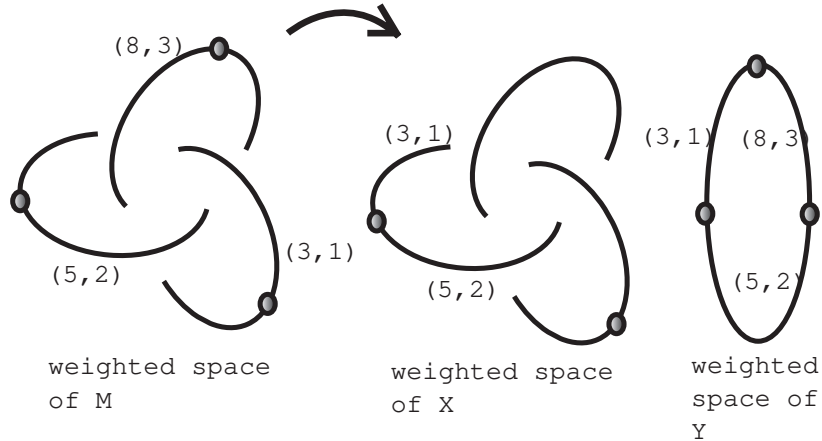}
\end{center}
\caption{An example of the decomposition of a weighted circle}
\label{fig3}
\end{figure}

Finally we consider the case that $E^\ast\cup F^\ast$ consists of
one weighted circle with three isolated fixed points lying on it. To
deal with this case, we decompose the weighted circle into two
weighted circles satisfying the following conditions (see
\cite{B02}):

\begin{enumerate}
\item[(1)]
$M$ is equivariantly diffeomorphic to the connected sum of $X$ and
$Y$. Here we do not require that $X$ and $Y$ admit a positively curved metric.

\item[(2)]
The weighted orbit space of $X$ is the same as $M$ except that
the weighted circle has exactly two isolated fixed points.

\item[(3)]
The weighted orbit space of $Y$ is $S^3$ with a trivially
embedded multiply-weighted circle with the original weights which is
unknotted in the orbit space $Y^\ast$ of $Y$. See Figure \ref{fig3}
\end{enumerate}

By considering the Euler characteristic of the fixed point set, we
can know that $X$ (resp. $Y$) is equivariantly homeomorphic to $S^4$
(resp. ${\bf CP}^2$). Then we again use the replacement trick of P.
Pao (Theorem \ref{thm3.1}) to obtain a new manifold $X'$
equivariantly diffeomorphic to $X$ either whose $E^\ast\cup F^\ast$
consists of two isolated fixed points or whose orbit space is the
orbit space $X^\ast$ minus the interior of 3-ball with $E^\ast\cup
F^\ast$ removed. But then by Theorem \ref{thm2.1} the $S^1$-action
of $X'$ lifts to $S^1\times S^1$-action. Hence $X'$ should be
equivariantly diffeomorphic to $S^4$.

Note that by construction $Y$ is equivariantly homeomorphic to ${\bf
CP}^2$ and has a weighted circle with three isolated fixed points
which is unknotted in $Y^\ast$ homeomorphic to $S^3$. Thus by
Theorem \ref{thm2.1} (Theorem 7.1 in \cite{F77}), the action of
$S^1$ on $Y$ extends to an action of $T^2=S^1\times S^1$. Now it is
a standard argument to show that $Y$ is equivariantly diffeomorphic
to ${\bf CP}^2$. To see it more precisely, note first that by
construction $Y$ acted on by $T^2$ does not have non-trivial finite
isotropy groups. (In fact, this is automatically true: if a closed
smooth orientable 4-manifold with an effective $T^2$-action whose
first homology group over integer coefficients is trivial then the
only finite isotropy group is the identity. See Lemma 5.2 in
\cite{OR70}.) Thus it follows from the classification of Orlik and
Raymond in \cite{OR70} that $Y$ should be $T^2$-equivariantly
homeomorphic to ${\bf CP}^2$. This in turn implies that their orbit
spaces are homeomorphic to each other in a weight-preserving way.
Since their orbit spaces are all 2-disks, they are indeed
weight-preservingly diffeomorphic. Thus by Theorem 4.2 in
\cite{OR70} we can conclude that there exists an equivariant
diffeomrphism of $Y$ onto ${\bf CP}^2$. This completes the proofs of
Theorems \ref{thm1.1} and \ref{thm1.2}.

\smallskip\smallskip
\noindent{\bf Acknowledgements:}
The author would like to thank the anonymous reader for providing valuable comments on this paper. This work was
supported by the Korea Research Foundation Grant and KOSEF.
\smallskip\smallskip

% ----------------------------------------------------------------

\bigskip

\address

\end{document}